\documentclass[12pt,leqno]{article}
\usepackage{amssymb,graphics,latexsym}
\usepackage{graphicx}
\usepackage{amsmath}
\large\normalsize

\def\R{\hbox{\bf\rlap{I}{\hbox to 2 pt{}}R}}

\newcommand{\re}{{\rm Re\, }}
\newcommand{\im}{{\rm Im\, }}

\newcommand{\dia}{{\rm diag\, }}

\begin{document}
\thispagestyle{empty}
\begin{center}
\section*{Higher-rank Numerical Ranges and Kippenhahn Polynomials}

\vspace*{8mm}

\begin {tabular}{lcl}
\hspace*{1cm}{\bf Hwa-Long Gau}$^1$\hspace*{1cm}&and & \hspace*{1cm}{\bf Pei Yuan Wu}$^2$
\vspace*{3mm}\\
Department of Mathematics & & Department of Applied
Mathematics\\
National Central University&& National Chiao Tung University\\
Chung-Li 320, Taiwan&& Hsinchu 300, Taiwan\\
hlgau@math.ncu.edu.tw&&pywu@math.nctu.edu.tw
\end{tabular}

\end{center}

\vspace{13mm}

\noindent
{\bf Abstract}

We prove that two $n$-by-$n$ matrices $A$ and $B$ have their rank-$k$ numerical ranges $\Lambda_k(A)$ and $\Lambda_k(B)$ equal to each other for all $k$, $1\le k\le \lfloor n/2\rfloor+1$, if and only if their Kippenhahn polynomials $p_A(x,y,z)\equiv\det(x\re A+y\im A+zI_n)$ and $p_B(x,y,z)\equiv\det(x\re B+y\im B+zI_n)$ coincide. The main tools for the proof are the Li-Sze characterization of higher-rank numerical ranges, Weyl's perturbation theorem for eigenvalues of Hermitian matrices and B\'{e}zout's theorem for the number of common zeros for two homogeneous polynomials.

\vspace{1cm}

\noindent
\emph{AMS Subject Classification}: 15A60\\
\emph{Keywords}: Higher-rank numerical range; Kippenhahn polynomial.

\vspace{1cm}

${}^1$Research supported in part by the National Science Council of the
Republic of China under project NSC 101-2115-M-008-006.

${}^2$Research supported in part by the National Science Council of the
Republic of China under project NSC 101-2115-M-009-004 and by the MOE-ATU.

\newpage

For an $n$-by-$n$ complex matrix $A$, its \emph{rank-$k$
numerical range} ($1\leq k\leq n$) is, by definition,
\[
\Lambda_k(A)=\{\lambda\in\mathbb{C}: X^*AX=\lambda I_k \mbox{ for some $n$-by-$k$ matrix $X$ with
$X^*X=I_k$}\}.
\]
Motivated by investigations in connection with the quantum error
correction, researchers started to study the higher-rank numerical ranges
in \cite{1}. The research was then pursued in a flurry of papers
\cite{2,3,4,5,6,7,14}. It is now known that $\Lambda_k(A)$, $1\leq k\leq
n$, is always convex \cite{4}, and, moreover, it consists of
those $\lambda$'s in $\mathbb{C}$ for which $\re
(e^{-i\theta}\lambda)\leq \lambda_k(\re (e^{-i\theta}A))$ for all
real $\theta$ \cite[Theorem 2.2]{5}. Here, and for our later discussions, we use $\re
X=(X+X^*)/2$ and $\im X=(X-X^*)/(2i)$ to denote the \emph{real}
and \emph{imaginary parts} of a finite matrix $X$, and, for an
$n$-by-$n$ Hermitian matrix $Y$, $\lambda_1(Y)\geq \cdots \geq
\lambda_n(Y)$ denote its (ordered) eigenvalues. Note that the
rank-one numerical range $\Lambda_1(A)$ coincides with the
classical numerical range $W(A)=\{\langle Ax,x \rangle:
x\in\mathbb{C}^n, \|x\|=1\}$ of $A$, where $\langle \cdot,\cdot
\rangle$ and $\|\cdot\|$ are the standard inner product and its
associated norm in $\mathbb{C}^n$.

\vspace{4mm}

The purpose of this paper is to determine when two matrices of
the same size have all their higher-rank numerical ranges equal to
each other. The following is the main theorem, which provides the
answer.

\vspace{.5cm}

{\bf Theorem 1.}  {\em The following conditions are equivalent for
$n$-by-$n$ matrices $A$ and $B$}:

(a) {\em $\Lambda_k(A)=\Lambda_k(B)$ for all $k$}, {\em $1\leq k\leq
\lfloor n/2 \rfloor +1$},

(b) {\em $\det(x\re A+y\im A+z I_n)=\det(x\re B+y\im B+z I_n)$
for all complex $x$}, {\em $y$ and $z$}, {\em and}

(c) {\em the eigenvalues of $\re (e^{-i\theta}A)$ and $\re
(e^{-i\theta}B)$ coincide} ({\em with the same
multiplicities}) {\em for all real $\theta$}.

\vspace{.5cm}

Here $\lfloor n/2 \rfloor$ denotes the largest integer which is
less than or equal to $n/2$. For an $n$-by-$n$ matrix $X$, we
call $p_X(x,y,z)=\det(x\re X+y\im X+z I_n)$ the \emph{Kippenhahn
polynomial} of $X$. It is a degree-$n$ homogeneous polynomial in
$x$, $y$ and $z$ with real coefficients.

\vspace{4mm}

Note that when the $n$-by-$n$ matrices $A$ and $B$ are such that
$p_A$ or $p_B$ is irreducible, the equality of $\Lambda_1(A)$ and
$\Lambda_1(B)$ already guarantees that $p_A$ and $p_B$ coincide
(cf. \cite[Corollary 2.4]{8}). On the other hand, the number
$\lfloor n/2 \rfloor +1$ in Theorem 1 (a) cannot be further
reduced as the 3-by-3
matrices $A=\dia (0,1,1)$ and $B=\dia (0,0,1)$ with
$\Lambda_1(A)=\Lambda_1(B)=[0,1]$, $p_A(x,y,z)=z(x+z)^2$ and
$p_B(x,y,z)=z^2(x+z)$ show. Also, the conditions in Theorem 1
cannot be strengthened to the unitary equivalence of $A$ and $B$. For example, if
\[
A=\begin{bmatrix}
0   & 1   &             \\
    & 0     & 2  \\
    &       & 0
\end{bmatrix}
\hspace{5mm}\mbox{and}\hspace{5mm}
B =
\begin{bmatrix}
0   & 2   &     \\
    & 0     & 1 \\
    &       & 0
\end{bmatrix},
\]
then $\Lambda_1(A)=\Lambda_1(B)=\{z\in\mathbb{C}: |z|\leq
\sqrt{5}/2\}$, $\Lambda_2(A)=\Lambda_2(B)=\{0\}$,
$\Lambda_3(A)=\Lambda_3(B)=\emptyset$ and
$p_A(x,y,z)=p_B(x,y,z)=z^3-(5/4)(x^2+y^2)z$, but $A$ and $B$ are
not unitarily equivalent (cf. \cite[Example 4]{9}). However, for
certain special-type matrices, we do have the unitary equivalence.

\vspace{.5cm}

{\bf Corollary 2.}  {\em Let $A$ and $B$ be $n$-by-$n$ matrices}. {\em If
$n=2$ or $A$ and $B$ are both normal or both companion matrices},
{\em then the conditions in Theorem {\em 1} are equivalent to the
unitary equivalence of $A$ and $B$}.

\vspace{.5cm}

{\em Proof}. It is well-known that if $X=\begin{bmatrix}
a    & c     \\
0    & b
\end{bmatrix}$, then $\Lambda_1(X)$ equals the elliptic disc with
foci $a$ and $b$ and with minor axis of length $|c|$. Hence, for
2-by-2 matrices, $\Lambda_1(A)=\Lambda_1(B)$ implies the unitary
equivalence of $A$ and $B$. For general matrices, if
$p_A(x,y,z)=p_B(x,y,z)$ for all $x$, $y$ and $z$, then plugging
in $x=1$ and $y=i$ yields $\det(A+zI_n)=\det(B+zI_n)$ for all
$z$, which implies that the eigenvalues of $A$ and $B$ coincide
(with the same algebraic multiplicities). In particular, if $A$
and $B$ are normal or companion matrices, then obviously they are
unitarily equivalent (in fact, equal in the latter case).
\hfill $\blacksquare$

\vspace{.5cm}

To prepare for the proof of Theorem 1, we review some basic properties of the Kippenhahn polynomials and numerical ranges. Recall that the \emph{complex projective plane} $\mathbb{CP}^2$ consists of the equivalence classes of ordered triple $[x,y,z]$ of complex numbers $x$, $y$ and $z$ which are not all equal to zero under the equivalence relation: $[x, y, z]\sim [x',y',z']$ if $[x,y,z]=\lambda[x',y',z']$ for some nonzero scalar $\lambda$. The point $[x,y,z]$ in $\mathbb{CP}^2$ with $z\neq 0$ corresponds to the point $(x/z, y/z)$ in $\mathbb{C}^2$ and, conversely, $(u,v)$ in $\mathbb{C}^2$ corresponds to $[u,v,1]$ in $\mathbb{CP}^2$. For a homogeneous polynomial $p$ in $x, y$ and $z$, the \emph{dual} of the algebraic curve $p(x,y,z)=0$ in $\mathbb{CP}^2$ is the curve
\[
\{[u,v,w]\in \mathbb{CP}^2: ux+vy+wz=0 \mbox{ is a tangent line of $p(x,y,z)=0$}\}.
\]
It is known that the dual of the dual curve is the original one (cf. \cite[Theorem 1.5.3]{10}). They have bearings on the numerical range because of Kippenhahn's result \cite{11}: \emph{the numerical range $W(A)$ of an $n$-by-$n$ matrix $A$ equals the convex hull of the real points $(u/w, v/w)$ of the dual curve of $p(x,y,z)=0$}. In the following, we also need B\'{e}zout's theorem \cite[Theorem 3.9]{12}, which counts the number of intersection points of two algebraic curves: \emph{if two homogeneous polynomials $p$ and $q$ in $x$, $y$ and $z$ of degrees $m$ and $n$}, \emph{respectively}, \emph{have no common factor}, \emph{then the number of common zeros of $p$ and $q$ is at most $mn$}.

\vspace{4mm}

We now proceed to prove Theorem 1. The main part is to show the implication (a)$\Rightarrow $(b). This is done via a series of lemmas. Note that the Kippenhahn polynomial $p_A$ of an $n$-by-$n$ matrix $A$ can be factored as the product of irreducible (real homogeneous) polynomials: $p_A=q_1^{n_1}\cdots q_m^{n_m}$, where the $q_j$'s are distinct and $n_j \geq 1$ for all $j$. Under the condition in Theorem 1 (a), we will show that each $q_j^{n_j}$ is also a factor of $p_B$ and thus $p_A$ divides $p_B$. Condition (b) then follows by symmetry. We start with the following lemma dealing with an irreducible $q$ having degree at least two.

\vspace{.5cm}

{\bf Lemma 3.} \emph{Let $q$ be an irreducible real homogeneous polynomial in $x, y$ and $z$ with degree at least two}. \emph{If $C$ is the curve in the plane consisting of the real points of the dual of $q(x,y,z)=0$}, \emph{then the convex hull of $C$ has no corner}.

\vspace{4mm}

Recall that, for a nonempty compact convex subset $\bigtriangleup$ of the plane, a point $\lambda$ on the boundary of $\bigtriangleup$ is a \emph{corner} of $\bigtriangleup$ if $\bigtriangleup$ has more than one supporting lines passing through it; otherwise, $\lambda$ is a \emph{differentiable point} of $\partial\bigtriangleup$.

\vspace{4mm}

{\em Proof of Lemma $3$}. Let $\lambda$ be a corner of the convex hull $\bigtriangleup$ of $C$. Then there are some $\theta_1$ and $\theta_2$, $\theta_1<\theta_2$, such that $x\cos\theta+y\sin\theta=\re(e^{-i\theta}\lambda)$ is a supporting line of $\bigtriangleup$ for all $\theta$ in $(\theta_1,\theta_2)$. By duality, this implies that $q(\cos\theta, \sin\theta, -\re(e^{-i\theta}\lambda))=0$ for all such $\theta$'s. On the other hand, $[\cos\theta, \sin\theta, -\re(e^{-i\theta}\lambda)]$ is also a zero of the linear polynomial $\lambda_1 x+\lambda_2y+z$, where $\lambda_1=(\lambda+\overline{\lambda})/2$ and $\lambda_1=(\lambda-\overline{\lambda})/(2i)$. B\'{e}zout's theorem then implies that $\lambda_1 x+\lambda_2y+z$ is a factor of $q$, which contradicts the irreducibility of $q$. Hence $\bigtriangleup$ cannot have any corner. \hfill $\blacksquare$

\vspace{4mm}

Essentially the same arguments as above were used in \cite{11} to prove that every corner of $W(A)$ for a finite matrix $A$ is an eigenvalue of $A$.

\vspace{4mm}

Another observation which we need is the following lemma, whose proof we omit.

\vspace{4mm}

{\bf Lemma 4.} \emph{Let $\bigtriangleup$ be a nonempty nonsingleton compact convex subset of the plane}. \emph{Let} $\bigtriangleup=\cap_{\theta\in[0,2\pi)}H_{\theta}$, \emph{where} $H_{\theta}=\{x+iy\in\mathbb{C}: x\cos\theta+y\sin\theta\le d(\theta)\}$ \emph{with} $\theta \mapsto d(\theta)$ \emph{continuous}, \emph{and let $\lambda$ be a point in the boundary of $\bigtriangleup$}.

(a) \emph{If $\lambda$ is a differentiable point of $\partial\bigtriangleup$}, \emph{then some $\partial H_{\theta}$ is the unique supporting line of $\bigtriangleup$ which passes through $\lambda$}.

(b) \emph{If $\bigtriangleup$ is not a line segment and $\lambda$ is a corner of $\bigtriangleup$}, \emph{then there are $\theta_1$ and $\theta_2$ in $[0,2\pi)$ with $\theta_1<\theta_2$ such that $\partial H_{\theta_1}$ and $\partial H_{\theta_2}$ are supporting lines of $\bigtriangleup$ and $\partial H_{\theta_1}\cap\partial H_{\theta_2}=\{\lambda\}$}. \emph{If we further require that} $(\theta_2-\theta_1)$ (\emph{mod} $\pi$) \emph{be maximal}, \emph{then $\theta_1$ and $\theta_2$ are unique}.

(c) \emph{If $\bigtriangleup$ is a line segment and $\lambda$ is an endpoint of $\bigtriangleup$}, \emph{then there are unique $\theta_1$ and $\theta_2$ in $[0,2\pi)$ with $\theta_2-\theta_1=\pi$ such that $\bigtriangleup\subseteq\partial H_{\theta_1}\cap\partial H_{\theta_2}$}.

\vspace{4mm}

In the following, this will be applied for an $n$-by-$n$ matrix $A$ with $\bigtriangleup=\Lambda_k(A)$ and $H_{\theta}=\{x+iy\in\mathbb{C}: x\cos\theta+y\sin\theta\le \lambda_k(\re(e^{-i\theta}A))\}$, $1\le k\le n$. Note that in this case $\partial H_{\theta}$ is in general not a supporting line of $\bigtriangleup$ nor the converse. One example is
$$A=\left[\begin{array}{cccc} -1/\sqrt{2} & -1/2 & 1/(2\sqrt{2}) & 1/4 \\ & -1/\sqrt{2} & -1/2 & -1/(2\sqrt{2}) \\ & & 1/\sqrt{2} & -1/2 \\ & & & 1/\sqrt{2} \end{array}\right]$$
as it is known that $\Lambda_2(A)$ has exactly two corners, around which the $\partial H_{\theta}$'s and the supporting lines are completely different (cf. \cite[Example 7]{15}).

\vspace{4mm}

We now proceed to obtain a characterization of $p_A$ with the power of an irreducible factor with degree at least two in terms of the relative positions of the $\Lambda_k(A)$'s.

\vspace{4mm}

{\bf Lemma 5.} \emph{Let $A$ be an $n$-by-$n$ matrix}, \emph{$q$ be an irreducible real homogeneous polynomial in $x, y$ and $z$ with degree at least two}, \emph{and $C$ be the real part of the dual curve of $q(x,y,z)=0$}. \emph{Then $q^m$ divides $p_A$ $(m\ge 1)$ if and only if $\partial\Lambda_{k_0}(A)\cap\partial\Lambda_{k_0-1}(A)\cap\cdots\cap\partial\Lambda_{k_0-m+1}(A)$ contains an arc of $C$ for some $k_0$}, $1\le k_0\le\lfloor n/2\rfloor$.

\vspace{4mm}

A result we need in the proof is Weyl's perturbation theorem for (ordered) eigenvalues of Hermitian matrices (cf. \cite[Theorem VI.2.1]{13}): \emph{if $X$ and $Y$ are $n$-by-$n$ Hermitian matrices with eigenvalues $\lambda_1(X)\ge\cdots\ge\lambda_n(X)$ and $\lambda_1(Y)\ge\cdots\ge\lambda_n(Y)$}, \emph{respectively}, \emph{then $|\lambda_j(X)-\lambda_j(Y)|\le\|X-Y\|$ for all $j$}, $1\le j\le n$.

\vspace{4mm}

{\em Proof of Lemma $5$}. Let $\bigtriangleup$ be the convex hull of $C$. Assume first that $q^m$ divides $p_A$. Let $k_0$ be the largest integer for which $\bigtriangleup$ is contained in $\Lambda_{k_0}(A)$. Since $\bigtriangleup\subseteq\Lambda_1(A)$ by Kippenhahn's result, we have $k_0\ge 1$. On the other hand, since $\Lambda_{\lfloor n/2\rfloor+1}(A)$ is either a singleton or an empty set \cite[Proposition 2.2]{1}, if $k_0>\lfloor n/2\rfloor$, then $\bigtriangleup\subseteq\Lambda_{\lfloor n/2\rfloor+1}(A)$ and hence $\bigtriangleup$ is a singleton. By duality, this says that $q$ is of degree one, contradicting our assumption that $q$ has degree at least two. Thus $1\le k_0\le\lfloor n/2\rfloor$. Note that, by Lemma 3, $\bigtriangleup$ has no corner. Hence, for each real $\theta$, $\bigtriangleup$ has a unique supporting line $x\cos\theta+y\sin\theta=d(\theta)$ with $x\cos\theta+y\sin\theta\le d(\theta)$ for all $x+iy$ in $\bigtriangleup$. Since $\bigtriangleup$ is not contained in $\Lambda_{k_0+1}(A)$, there is some $\lambda_0$ in $\bigtriangleup\setminus\Lambda_{k_0+1}(A)$. By the Li-Sze characterization \cite[Theorem 2.2]{5} of $\Lambda_{k_0+1}(A)$, we have $\re(e^{-i\theta_0}\lambda_0)>\alpha_{k_0+1}(\theta_0)$ for some $\theta_0$. Here $\alpha_k(\theta)$ denotes $\lambda_k(\re(e^{-i\theta}A))$ for $1\le k\le n$ and $\theta$ in $\mathbb{R}$. Weyl's  perturbation theorem then implies that there is some $\delta>0$ such that $\re(e^{-i\theta}\lambda_0)>\alpha_{k_0+1}(\theta)$ for all $\theta$ in $(\theta_0-\delta, \theta_0+\delta)$. On the other hand, since $\lambda_0$ is in $\bigtriangleup$ and $\bigtriangleup$ is contained in $\Lambda_{k_0}(A)$, we also have $\re(e^{-i\theta}\lambda_0)\le d(\theta)\le\alpha_{k_0}(\theta)$ for all $\theta$. Thus $\alpha_{k_0+1}(\theta)<d(\theta)\le\alpha_{k_0}(\theta)$ for all $\theta$ in $(\theta_0-\delta, \theta_0+\delta)$. Since $[\cos\theta,\sin\theta,-d(\theta)]$ is a zero of $q(x,y,z)$ by duality, the fact that $q^m$ divides $p_A$ implies that $[\cos\theta,\sin\theta,-d(\theta)]$ is a zero of $p_A(x,y,z)$ with multiplicity at least $m$. We infer from above that $d(\theta)=\alpha_k(\theta)$ for all $k$, $k_0-m+1\le k\le k_0$, and all $\theta$ in $(\theta_0-\delta, \theta_0+\delta)$. We then obtain from $\bigtriangleup\subseteq\Lambda_{k}(A)$ that $x\cos\theta+y\sin\theta=d(\theta)$ is the unique supporting line of $\Lambda_k(A)$ for all such $k$'s and $\theta$'s. Hence $\partial\Lambda_{k_0}(A)\cap\partial\Lambda_{k_0-1}(A)\cap\cdots\cap\partial\Lambda_{k_0-m+1}(A)$ contains an arc of $C$.

\vspace{4mm}

For the converse, if $\partial\Lambda_{k_0}(A)\cap\cdots\cap\partial\Lambda_{k_0-m+1}(A)$ contains an arc of $C$, then the supporting line $x\cos\theta+y\sin\theta=d(\theta)$ of $\bigtriangleup$ is also a supporting line of $\Lambda_k(A)$ for all $k$, $k_0-m+1\le k\le k_0$, and all $\theta$ in some $(\theta_1, \theta_2)$. This implies by the Li-Sze characterization \cite[Theorem 2.2]{5} of $\Lambda_{k}(A)$ and Lemma 4 (a) that $d(\theta)=\alpha_k(\theta)$. Hence $[\cos\theta,\sin\theta,-d(\theta)]$ is a zero of $p_A(x,y,z)$ with multiplicity at least $m$ for all such $\theta$'s. Since $[\cos\theta,\sin\theta,-d(\theta)]$ is a zero of $q(x,y,z)$ for all $\theta$ by duality, we obtain that $p_A(x,y,z)$ and $q(x,y,z)$ have infinitely many common zeros of the form $[\cos\theta,\sin\theta,-d(\theta)]$. B\'{e}zout's theorem yields that the irreducible $q$ divides $p_A$. Next we claim that the number of $\theta$'s in $[0,2\pi)$ for which $-d(\theta)$ is a zero of $q(\cos\theta,\sin\theta,z)$ with multiplicity at least two is finite. Indeed, if otherwise, then $p(x,y,z)\equiv\partial q(x,y,z)/\partial z$ and $q(x,y,z)$ have infinitely many common zeros of the form $[\cos\theta,\sin\theta,-d(\theta)]$. Since $q$ is irreducible, B\'{e}zout's theorem implies that $q$ divides $p$, which is impossible for the degree of $q(\cos\theta,\sin\theta,z)$ is one bigger than that of $p(\cos\theta,\sin\theta,z)$. Hence again we can apply B\'{e}zout's theorem to $p_A/q$ and $q$ to obtain that $q$ divides $p_A/q$. Repeating these arguments yield that $q$ divides $p_A/q^j$ for all $j$, $0\le j\le m-1$. Thus $q^m$ divides $p_A$, completing the proof. \hfill $\blacksquare$

\vspace{4mm}

{\bf Corollary 6.} \emph{An $n$-by-$n$ matrix $A$ is normal if and only if $\Lambda_k(A)$ is a} (\emph{closed}) \emph{polygonal region for all $k$}, $1\le k\le n$.

\vspace{4mm}

Here a polygonal region is one whose boundary is a polygon. In the degenerate case, this may be an empty set, a singleton or a line segment.

\vspace{4mm}

{\em Proof of Corollary $6$}. The necessity follows from \cite[Corollary 2.4]{5}. The sufficiency is an easy consequence of Lemma 5 since the latter implies that $p_A$ has only linear factors. \hfill $\blacksquare$

\vspace{4mm}

The next corollary is another consequence of Lemma 5.

\vspace{4mm}

{\bf Corollary 7.} {\em Let $A$ and $B$ be $n$-by-$n$ matrices with $\Lambda_k(A)=\Lambda_k(B)$ for all $k$}, {\em $1\leq k\leq
\lfloor n/2 \rfloor +1$}. {\em Then $p_A$ and $p_B$ contain the same powers of irreducible factors with degrees at least two}.

\vspace{4mm}

To show that $p_A$ and $p_B$ contain the same powers of linear factors under the above conditions, we need a characterization, analogous to the one in Lemma 5, for powers of linear factors. Unfortunately, a complete analogue of Lemma 5 is not true. We have only had the following necessary condition.

\vspace{4mm}

{\bf Lemma 8.} {\em Let $A$ be an $n$-by-$n$ matrix}. {\em If $(ax+by+z)^m$ divides $p_A(x,y,z)$}, {\em where $a$ and $b$ are real and $m\geq 1$}, {\em then there is a $k_0$}, {\em $k_0\geq m$}, {\em such that $a+bi$ is a corner of $\Lambda_k(A)$ for all $k$}, {\em $k_0-m+1\leq k\leq k_0$}, {\em and is not in $\Lambda_{k_0+1}(A)$}.

\vspace{4mm}

{\em Proof}. Let $k_0=\max\{k\geq 1: a+bi \mbox{ is in } \Lambda_k(A)\}$. Then $a+bi$ is in $\Lambda_{k}(A)$ for all $k$, $1\le k\le k_0$, and is not in $\Lambda_{k_0+1}(A)$. Because $a+bi$ is not in $\Lambda_{k_0+1}(A)$, there is some $\theta_0$ such that $a\cos\theta_0+b\sin\theta_0>\alpha_{k_0+1}(\theta_0)\equiv \lambda_{k_0+1}(\re (e^{-i\theta_0}A))$. By Weyl's perturbation theorem, we obtain $a\cos\theta+b\sin\theta>\alpha_{k_0+1}(\theta)$ on $(\theta_0-\delta, \theta_0+\delta)$ for some $\delta>0$, where, for $1\leq k\leq n$ and $\theta$ in $\mathbb{R}$, $\alpha_k(\theta)$ denotes $\lambda_k(\re (e^{-i\theta}A))$. On the other hand, $a+bi$ being in $\Lambda_{k_0}(A)$ implies that $a\cos\theta+b\sin\theta\leq\alpha_{k_0}(\theta)$ for all real $\theta$. Since $[\cos\theta, \sin\theta,-(a\cos\theta+b\sin\theta)]$ is a zero of $(ax+by+z)^m$ and hence a zero of $p_A(x,y,z)$ with multiplicity at least $m$, $a\cos\theta+b\sin\theta$ is an eigenvalue of $\re (e^{-i\theta}A)$ with multiplicity at least $m$. Hence $k_0\ge m$ and $a\cos\theta+b\sin\theta=\alpha_k(\theta)$ for all $k$, $k_0-m+1\leq k\leq k_0$, on $(\theta_0-\delta, \theta_0+\delta)$. This means that $a+bi$ is a corner of $\Lambda_k(A)$ for all such $k$'s.  \hfill $\blacksquare$

\vspace{4mm}

The next example shows that the converse of the assertion in Lemma 8 is not necessarily true. Recall that, for any subset $\bigtriangleup$ of the plane, $\bigtriangleup^{\wedge}$ denotes its convex hull.

\vspace{4mm}

{\bf Example 9.} Let $A=\dia (1, i, -1, -i, 1/2, i/2, -1/2, -i/2, (1+i)/3)$. Then $\Lambda_1(A)=\{\pm 1,\pm i\}^{\wedge}$, $\Lambda_2(A)=\{\pm 1/2,
\pm i/2, (\pm 1\pm i)/3\}^{\wedge}$, $\Lambda_3(A)=\{0, 1/4, i/4, (1+i)/3\}^{\wedge}$, $\Lambda_4(A)=\{0\}$ and $\Lambda_k(A)=\emptyset$ for $5\leq k\leq 9$ (cf. Figure 10). Hence $(1+i)/3$ is a corner of $\Lambda_2(A)$ and $\Lambda_3(A)$, but $((1/3)x+(1/3)y+z)^2$ does not divide $p_A(x,y,z)=(x^2-z^2)(y^2-z^2)((1/4)x^2-z^2)((1/4)y^2-z^2)((1/3)x+(1/3)y+z)$.

\vspace{2mm}

\centerline{\includegraphics[scale=0.56]{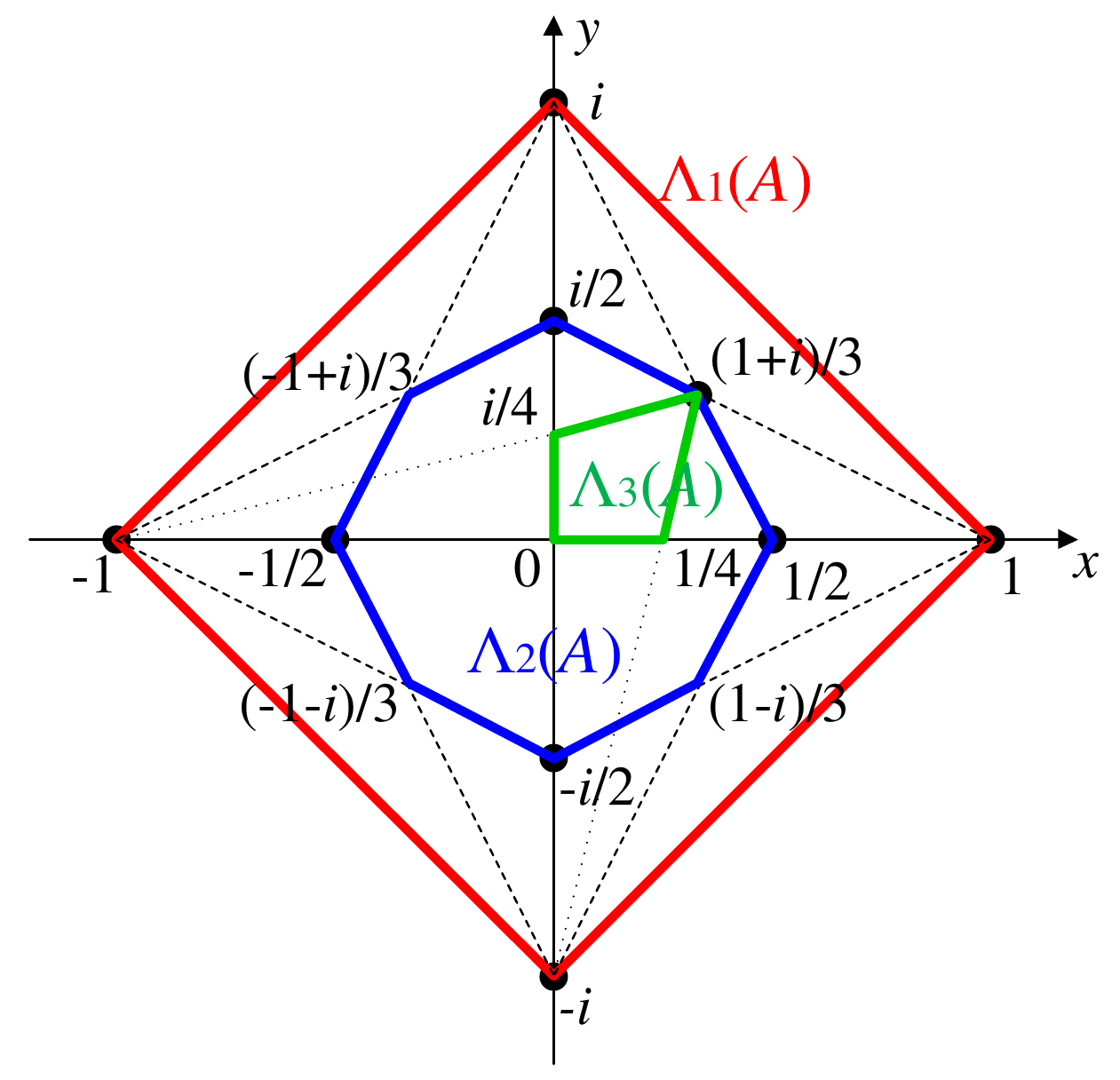}}
\centerline{Figure 10}

\vspace{4mm}

For an $n$-by-$n$ matrix $A$ and $1\le\ell\le\lfloor n/2 \rfloor$, let
\begin{eqnarray*}
V_{\ell}(A)&=&\{a+bi : ax+by+z \mbox{ is a real linear factor of } p_A(x,y,z) \mbox{ with multiplicity} \\
&& \ \ m, 1\le m\le\lfloor n/2 \rfloor-\ell +1, \mbox{ and } a+bi\in\Lambda_{\ell+m-1}(A)\setminus\Lambda_{\ell+m}(A)\}.
\end{eqnarray*}
We remark that if $ax+by+z$ is a real linear factor of $p_A(x,y,z)$ with multiplicity $m$, then Lemma 8 yields that there is a $k_0$, $k_0\geq m$, such that $a+bi$ is a corner of $\Lambda_k(A)$ for all $k$, $k_0-m+1\leq k\leq k_0$, and is not in $\Lambda_{k_0+1}(A)$. If, moreover, $k_0\le \lfloor n/2\rfloor$ or, equivalently, $a+bi\not\in\Lambda_{\lfloor n/2 \rfloor +1}(A)$,  then $a+bi$ is in $V_{k_0-m+1}(A)$. Conversely, if $a+bi$ is in $V_{\ell_0}(A)$ and $a+bi\in\Lambda_{k_0}(A)\setminus\Lambda_{k_0+1}(A)$, then the definition of $V_{\ell_0}(A)$ yields that $a+bi\not\in\Lambda_{\lfloor n/2 \rfloor +1}(A)$, $\ell_0\le k_0\le\lfloor n/2 \rfloor$ and $ax+by+z$ is a real linear factor of $p_A(x,y,z)$ with multiplicity $k_0-\ell_0+1$. Obviously, the $V_{\ell}(A)$'s and $\Lambda_{\lfloor n/2 \rfloor +1}(A)$ are mutually disjoint and $V_1(A)\cup V_2(A)\cup\cdots\cup V_{\lfloor n/2 \rfloor}(A)\cup\Lambda_{\lfloor n/2 \rfloor +1}(A)=\{a+bi : ax+by+z \mbox{ is a real linear factor of } p_A(x,y,z)\}$.
For the proof of the latter, note that $\Lambda_{\lfloor n/2 \rfloor +1}(A)$ is either empty or a singleton (cf. \cite[Proposition 2.2]{1}). Hence we need only show that if $\Lambda_{\lfloor n/2 \rfloor +1}(A)=\{a+bi\}$, then $ax+by+z$ is a factor of $p_A(x,y,z)$. Indeed, in this case, $[\cos\theta, \sin\theta, -\lambda_{\lfloor n/2 \rfloor +1}(\re(e^{-i\theta}A))]$ is a zero of both $p_A(x,y,z)$ and $ax+by+z$ for all real $\theta$ by the Li-Sze characterization of $\Lambda_{\lfloor n/2 \rfloor +1}(A)$. The assertion then follows from B\'{e}zout's theorem.
As an example, in Example 9, we have $V_1(A)=\{\pm 1, \pm i\}$, $V_2(A)=\{\pm 1/2, \pm i/2\}$, $V_3(A)=\{(1+i)/3\}$ and $V_{4}(A)=\emptyset$.

\vspace{4mm}

The next lemma would help us to conclude the proof of Theorem 1.

\vspace{4mm}

{\bf Lemma 11.} {\em Let $A$ and $B$ be $n$-by-$n$ matrices with $\Lambda_k(A)=\Lambda_k(B)$ for all $k$}, {\em $1\leq k\leq
\lfloor n/2 \rfloor +1$}. {\em Then $V_{\ell}(A)=V_{\ell}(B)$ and $\lambda_{\ell}(\re(e^{-i\theta}A))=\lambda_{\ell}(\re(e^{-i\theta}B))$ for all $\ell$}, $1\le \ell\le\lfloor n/2 \rfloor$, and $\theta$ in $\mathbb{R}$.

\vspace{4mm}

{\em Proof}. For $1\le k\le n$ and real $\theta$, let $\alpha_{k}(\theta)=\lambda_{k}(\re(e^{-i\theta}A))$ and $\beta_{k}(\theta)=\lambda_{k}(\re(e^{-i\theta}B))$. We now prove our assertion by induction on $\ell$.

\vspace{3mm}

If $\ell=1$, then, since $A$ and $B$ have the same numerical range, we have $\alpha_1(\theta)=\beta_1(\theta)$ for all real $\theta$. We now check that $V_1(A)\subseteq V_1(B)$. Indeed, if $a+bi\in V_1(A)$ and $ax+by+z$ is a real linear factor of $p_A(x,y,z)$ with multiplicity $m$, $1\le m\le\lfloor n/2 \rfloor$, then Lemma 8 says that $a+bi$ is a corner of $\Lambda_k(A)$, $1\le k\le m$, and is not in $\Lambda_{m+1}(A)$. Since $\Lambda_k(A)=\Lambda_k(B)$ for all $k$, $1\leq k\leq \lfloor n/2 \rfloor +1$, and $m\le\lfloor n/2 \rfloor$, we obtain that $a+bi$ is also a corner of $\Lambda_k(B)$ for $1\leq k\leq m$ and is not in $\Lambda_{m+1}(B)$. As $a+bi$ is a corner of $\Lambda_1(B)$, it is an eigenvalue of $B$ and $ax+by+z$ divides $p_B(x,y,z)$. Moreover, there are $\theta_1$ and $\theta_2$, $\theta_1<\theta_2$, such that
$a\cos\theta+b\sin\theta=\beta_1(\theta)$ for all $\theta$ in $(\theta_1,\theta_2)$. On the other hand, $a+bi$ being in $\Lambda_m(B)$ implies that $a\cos\theta+b\sin\theta\leq\beta_m(\theta)$ for all real $\theta$ by the Li-Sze characterization of the $\Lambda_m(B)$. We thus have $a\cos\theta+b\sin\theta=\beta_1(\theta)=\cdots=\beta_m(\theta)$ and therefore $(a\cos\theta+b\sin\theta+z)^m$ divides $p_B(\cos\theta, \sin\theta, z)$ for $\theta$ in $(\theta_1,\theta_2)$. Using B\'{e}zout's theorem repeatedly, we conclude that $(ax+by+z)^m$ divides $p_B(x, y, z)$. This means that $ax+by+z$ is a real linear factor of $p_B(x,y,z)$ with multiplicity at least $m$. But since $a+bi\in\Lambda_m(B)\setminus\Lambda_{m+1}(B)$, $1\le m\le\lfloor n/2 \rfloor$, by Lemma 8 and the definition of $V_{1}(B)$, we deduce that $ax+by+z$ is a real linear factor of $p_B(x,y,z)$ with multiplicity $m$ and $a+bi$ is in $V_1(B)$. Therefore, we have $V_1(A)\subseteq V_1(B)$. By symmetry, we also obtain $V_1(B)\subseteq V_1(A)$. Hence we conclude that $V_1(A)= V_1(B)$.

\vspace{3mm}

Next assume that our assertion is true for all $\ell$, $1\le\ell<\ell_0$, and $\ell_0\le\lfloor n/2\rfloor$. We prove its validity for $\ell_0$. Firstly, we check that $\alpha_{\ell_0}(\theta)=\beta_{\ell_0}(\theta)$ for all real $\theta$. Indeed, if otherwise, then we have $\alpha_{\ell_0}(\theta_0)\neq\beta_{\ell_0}(\theta_0)$ for some $\theta_0$. Without loss of generality, we may assume that $\alpha_{\ell_0}(\theta_0)<\beta_{\ell_0}(\theta_0)$. Weyl's perturbation theorem then yields that $\alpha_{\ell_0}(\theta)< \beta_{\ell_0}(\theta)$ for all $\theta$ in some neighborhood $I\equiv (\theta_0-\delta, \theta_0+\delta)$ of $\theta_0$ ($\delta>0$). Since $\beta_{\ell_0}(\theta)$ is an eigenvalue of $\re (e^{-i\theta}B)$, we have $p_B(\cos\theta, \sin\theta, -\beta_{\ell_0}(\theta))=0$ for all real $\theta$. Note that $p_B$ has only finitely many irreducible factors. Thus there is some irreducible factor $q$ of $p_B$ and an infinite subset $I_1$ of $I$ such that $q(\cos\theta, \sin\theta, -\beta_{\ell_0}(\theta))=0$ for all $\theta$ in $I_1$.

\vspace{3mm}

If $q$ is of degree at least two, then, by Corollary 7, $q$ is also a factor of $p_A$. Hence $p_A(\cos\theta, \sin\theta, -\beta_{\ell_0}(\theta))=0$ and thus $\beta_{\ell_0}(\theta)$ is an eigenvalue of $\re (e^{-i\theta}A)$ for all $\theta\in I_1$. Moreover, by the induction hypothesis, we have $\alpha_{\ell_0-1}(\theta)=\beta_{\ell_0-1}(\theta)\ge\beta_{\ell_0}(\theta)>\alpha_{\ell_0}(\theta)$ for all $\theta\in I$, and thus $\beta_{\ell_0}(\theta)=\alpha_{\ell_0-1}(\theta)=\beta_{\ell_0-1}(\theta)>\alpha_{\ell_0}(\theta)$ and $q(\cos\theta, \sin\theta, -\alpha_{\ell_0-1}(\theta))=0$ for all $\theta$ in $I_1$. Let $m$ be the multiplicity of $q(x,y,z)$ in $p_A(x,y,z)$. Then $\alpha_{\ell_0}(\theta)<\alpha_{\ell_0-1}(\theta)=\cdots=\alpha_{\ell_0-m}(\theta)$ for all $\theta\in I_1$. The induction hypothesis says that $\alpha_{\ell}(\theta)=\beta_{\ell}(\theta)$ for all $\ell$, $1\le\ell<\ell_0$, and all real $\theta$, and thus $\beta_{\ell_0-m}(\theta)=\cdots=\beta_{\ell_0-1}(\theta)=\beta_{\ell_0}(\theta)$ for all $\theta\in I_1$. Since the set $I_1$ is infinite, using B\'{e}zout's theorem repeatedly (as in the proof of Lemma 5), we obtain that $q^{m+1}$ divides $p_B$. This means that $q$ is an irreducible factor of $p_B$ with multiplicity at least $m+1$, which contradicts the assertion of Corollary 7. Hence $q$ must be linear.

\vspace{3mm}

Let $q(x,y,z)=ax+by+z$. Then $\alpha_{\ell_0}(\theta)<\beta_{\ell_0}(\theta)=a\cos\theta+b\sin\theta$ for all $\theta\in I_1$, and hence $a+bi$ is not in $\Lambda_{\ell_0}(A)=\Lambda_{\ell_0}(B)$ by the Li-Sze characterization of $\Lambda_{\ell_0}(A)$. This implies that $a+bi$ is in $V_{\ell}(B)$ for some $\ell<\ell_0$. The induction hypothesis yields that $a+bi$ is also in $V_{\ell}(A)\ (=V_{\ell}(B))$. In particular, $ax+by+z$ is a real linear factor of $p_A(x,y,z)$. Hence $p_A(\cos\theta, \sin\theta, -\beta_{\ell_0}(\theta))=0$ and thus $\beta_{\ell_0}(\theta)$ is an eigenvalue of $\re (e^{-i\theta}A)$ for all $\theta\in I_1$. Moreover, by the induction hypothesis, we have $\alpha_{\ell_0-1}(\theta)=\beta_{\ell_0-1}(\theta)\ge\beta_{\ell_0}(\theta)>\alpha_{\ell_0}(\theta)$ for all $\theta\in I$. This forces that $\alpha_{\ell_0-1}(\theta)=\beta_{\ell_0}(\theta)>\alpha_{\ell_0}(\theta)$ for all $\theta\in I_1$. Let $m$ be the multiplicity of $ax+by+z$ in $p_A(x,y,z)$. Then $\alpha_{\ell_0}(\theta)<\alpha_{\ell_0-1}(\theta)=\cdots=\alpha_{\ell_0-m}(\theta)$ for all $\theta\in I_1$. The induction hypothesis implies that $\beta_{\ell_0-m}(\theta)=\cdots=\beta_{\ell_0-1}(\theta)=\beta_{\ell_0}(\theta)$ for all $\theta\in I_1$. Using B\'{e}zout's theorem repeatedly, we obtain that $(ax+by+z)^{m+1}$ divides $p_B(x,y,z)$. This means that $ax+by+z$ is a linear factor of $p_B(x,y,z)$ with multiplicity at least $m+1$. On the other hand, since $a+bi\in V_{\ell}(A)=V_{\ell}(B)$, under Lemma 8 and our assumption, $a+bi$ is a corner of $\Lambda_k(A)=\Lambda_k(B)$ for $\ell\leq k\leq \ell+m-1$ and is not in $\Lambda_{\ell+m}(A)=\Lambda_{\ell+m}(B)$. Hence the definition of $V_{\ell}(B)$ yields that the multiplicity of $ax+by+z$ in $p_B(x,y,z)$ is $m$, a contradiction. Thus $\alpha_{\ell_0}(\theta)= \beta_{\ell_0}(\theta)$ for all real $\theta$ as asserted.

\vspace{3mm}

We now show that $V_{\ell_0}(A)\subseteq V_{\ell_0}(B)$.
Suppose that $a+bi\in V_{\ell_0}(A)$ and $ax+by+z$ is a real linear factor of $p_A(x,y,z)$ with multiplicity $m$.
Lemma 8 yields that $a+bi$ is a corner of $\Lambda_k(A)$ for $\ell_0\leq k\leq k_0$ and is not in $\Lambda_{k_0+1}(A)$, where $k_0=\ell_0+m-1$. Note that $a+bi\not\in\Lambda_{\lfloor n/2\rfloor+1}(A)$ implies $k_0\le\lfloor n/2\rfloor$. Since $\Lambda_k(A)=\Lambda_k(B)$ for all $k$, $1\le k\le\lfloor n/2\rfloor+1$, $a+bi$ is also a corner of $\Lambda_k(B)=\Lambda_k(A)$ for $\ell_0\leq k\leq k_0$ and is not in $\Lambda_{k_0+1}(B)=\Lambda_{k_0+1}(A)$. The former implies that $a\cos\theta+b\sin\theta\leq \alpha_{k_0}(\theta), \beta_{k_0}(\theta)$ for all real $\theta$ while the latter, by Lemma 4 (b) and (c), that $a\cos\theta_0+b\sin\theta_0>\alpha_{k_0+1}(\theta_0), \beta_{k_0+1}(\theta_0)$ for some common $\theta_0$, both by the Li-Sze characterization of the higher-rank numerical ranges. Weyl's perturbation theorem then yields that $a\cos\theta+b\sin\theta>\alpha_{k_0+1}(\theta), \beta_{k_0+1}(\theta)$ for all $\theta$ in some neighborhood $I\equiv (\theta_0-\delta, \theta_0+\delta)$ of $\theta_0$ ($\delta>0$). Since $(ax+by+z)^m$ divides $p_A(x, y, z)$, $a\cos\theta+b\sin\theta$ appears as $m \, (=k_0-\ell_0+1)$ values of the $\alpha_k(\theta)$'s. Thus $a\cos\theta+b\sin\theta=\alpha_{k_0}(\theta)=\cdots=\alpha_{\ell_0}(\theta)$ for $\theta$ in $I$. Since we have proved that $\alpha_{\ell_0}(\theta)=\beta_{\ell_0}(\theta)$ for all real $\theta$, thus $a\cos\theta+b\sin\theta=\beta_{\ell_0}(\theta)$ for all $\theta$ in $I$. It follows that $a\cos\theta+b\sin\theta=\beta_{k_0}(\theta)=\cdots=\beta_{\ell_0}(\theta)$ for $\theta$ in $I$. Using B\'{e}zout's theorem repeatedly, we obtain that $(ax+by+z)^{m}$ divides $p_B$. This means that $ax+by+z$ is a real linear factor of $p_B(x,y,z)$ with multiplicity at least $m$. Moreover, from the definition of the $V_{\ell}(B)$'s, $a+bi\in\Lambda_{\ell_0+m-1}(B)\setminus\Lambda_{\ell_0+m}(B)$ implies that $a+bi\in V_{\ell_1}(B)$ for some $\ell_1\le\ell_0$. If $\ell_1<\ell_0$, the induction hypothesis yields that $a+bi$ is also in $V_{\ell_1}(A)$, which contradicts the mutual disjointness of the $V_{\ell}(A)$'s. Hence we conclude that $\ell_1=\ell_0$ or $a+bi\in V_{\ell_0}(B)$ as desired. For the converse, interchanging $A$ with $B$ in the above arguments, we also obtain $V_{\ell_0}(B)\subseteq V_{\ell_0}(A)$. Thus $V_{\ell_0}(A)= V_{\ell_0}(B)$, completing the proof.
 \hfill $\blacksquare$

\vspace{4mm}

We are now ready to prove Theorem 1.

\vspace{4mm}

{\em Proof of Theorem $1$}. The implication (b)$\Rightarrow$(c) is trivial and the implication (c)$\Rightarrow$(a) follows from the Li-Sze characterization of the higher-rank numerical ranges. We need only prove (a)$\Rightarrow $(b). Suppose that $q(x,y,z)$ is an irreducible factor of $p_A(x,y,z)$ with multiplicity $m$. If $q$ is of degree at least two, then Corollary 7 implies that $q(x,y,z)$ is also an irreducible factor of $p_B(x,y,z)$ with multiplicity $m$. Assume next that $q$ is of degree one, say, $q(x,y,z)=ax+by+z$. If $a+bi$ is not in $\Lambda_{\lfloor n/2 \rfloor +1}(A)$, then $a+bi$ is in $V_{\ell_0}(A)$ and $a+bi\in\Lambda_{\ell_0+m-1}(A)\setminus\Lambda_{\ell_0+m}(A)$ for some $\ell_0$ with $\ell_0+m\le\lfloor n/2 \rfloor+1$. Lemma 11 and the condition in (a) implies that $a+bi$ is also in $V_{\ell_0}(B)$ and $a+bi\in\Lambda_{\ell_0+m-1}(B)\setminus\Lambda_{\ell_0+m}(B)$. Hence $ax+by+z$ is also a factor of $p_B(x,y,z)$ with multiplicity $m$. Therefore, if $\Lambda_{\lfloor n/2 \rfloor +1}(A)$ is empty, then, since $p_A$ and $p_B$ have the same degree, we have $p_A=p_B$.

\vspace{3mm}

On the other hand, if $\Lambda_{\lfloor n/2 \rfloor +1}(A)$ is nonempty, then it must be a singleton, say, $\Lambda_{\lfloor n/2 \rfloor +1}(A)=\{c+di\}$. We have $\Lambda_{\lfloor n/2 \rfloor +1}(B)=\Lambda_{\lfloor n/2 \rfloor +1}(A)=\{c+di\}$, and $cx+dy+z$ is a real linear factor of both $p_A(x,y,z)$ and $p_B(x,y,z)$ (cf. the paragraph after Example 9). Since the degrees of $p_A$ and $p_B$ coincide, we infer from what were proved before that $p_A(x,y,z)=p_B(x,y,z)$, completing the proof.\hfill $\blacksquare$

\end{document}